\documentstyle[amscd]{amsart}
\numberwithin{equation}{section}
\theoremstyle{plain}
\newtheorem{thm}{Theorem}[section]
\newtheorem{cor}[thm]{Corollary}
\newtheorem{lem}[thm]{Lemma}
\newtheorem{prop}[thm]{Proposition}
\begin{document}
\title{A notion of synchronization of symbolic dynamics and a class of $C^*$-algebras}
\author{Wolfgang Krieger}
\address{Institute for Applied  Mathematics,
University of Heidelberg,
Im Neuenheimer Feld 294,
69120 Heidelberg, Germany}
\author{Kengo Matsumoto}
\address{ 
Department of Mathematics, 
Joetsu  University of Education,
Joetsu 943-8512 Japan}
\dedicatory{\large{\it{Dedicated to the memory of Ki Hang Kim}}}
\maketitle
\begin{abstract}
We discuss a synchronization property for subshifts,
that we call $\lambda$-synchronization.
Under an irreducibility assumption
we associate to a $\lambda$-synchronizing subshift
a simple and purely infinite $C^*$-algebra.
\end{abstract}

\def\Zp{{ {\Bbb Z}_+ }}
\def\M{{ {\cal M} }}
\def\LCHDN{{{{\frak L}^{Ch(D_N)}}}}
\def\Ext{{{\operatorname{Ext}}}}
\def\id{{{\operatorname{id}}}}
\def\Ker{{{\operatorname{Ker}}}}
\def\Card{{{\operatorname{Card}}}}
\def\card{{{\operatorname{card}}}}



\bigskip

\section{Introduction}
Let $\Sigma$ be a finite alphabet, and let $S_\Sigma $ be the left shift on
$\Sigma^{\Bbb Z}$,
$$
  S_\Sigma((x_i)_{ i \in {\Bbb Z}})   = (x_{i+1})_{ i \in {\Bbb Z}}, \qquad  
(x_i)_{ i \in {\Bbb Z}} \in \Sigma^{\Bbb Z}.
$$
The closed shift-invariant subsystems of the shift
 $(\Sigma^{\Bbb Z}, S_\Sigma) $ 
are called
subshifts. 
For an introduction to their theory, which belongs to symbolic
dynamics, we refer to \cite{Ki} and \cite{LM}. 
A finite word in the symbols of 
$\Sigma$  is called admissible for the subshift 
$X \subset \Sigma^{\Bbb Z} $ if it
appears somewhere in a point of $X$. 
A subshift is uniquely determined by its
language of admissible words that we denote by ${\cal L}(X)$. 
We let ${\cal L}_n(X)$ 
denote the set of words in ${\cal L}(X)$
of length $n \in {\Bbb N}$.

We set for a subshift $X \subset \Sigma^{\Bbb Z}$
\begin{equation*}
\Gamma^+(a) = \{ c \in {\cal L}(X) \mid a c  \in {\cal L}(X) \}, 
\qquad
a \in {\cal L}(X),
\end{equation*}
and
\begin{equation*}
\omega^+_n(a) = \bigcap_{c\in \Gamma^-(a)} \{ b \in {\cal L}_n(X) \mid c a b \in {\cal L}(X) \}, 
\qquad
a \in {\cal L}(X).
\end{equation*}
$\Gamma^-$ and $\omega^-$ have the time symmetric meaning.
An admissible word $v$ of a subshift $X \subset \Sigma^{\Bbb Z}$
is called a synchronizing word of $X$ 
if for $u,w \in {\cal L}(X)$ 
such that 
$uv, \ vw \in {\cal L}(X)$ 
also
$uvw \in {\cal L}(X)$. 
A topological transitive subshift
is said to be synchronizing if it has a synchronizing word.

In \cite{Kr6}
a property (D) of subshifts was introduced that expresses
a quality of synchronization 
that is weaker than synchronization
(For other notions of synchronization see  \cite{Kr2000},\cite{Kr8}, \cite{KM}).
A subshift $X\subset {\Sigma}^{\Bbb Z}$ 
has property (D)
if for  $\sigma \in \Sigma$ and 
$b \in \Gamma^-(\sigma)$ 
there exists a word $a \in \Gamma^-(b)$ such that
$\sigma \in \omega_1^+(ab)$.

Whereas synchronization of subshifts is time symmetric, 
property (D) is not. 
There exist coded systems \cite{BH} with property (D) whose inverse does  not  
have property (D) (see Section 4.3). 
Because of  the occurrence of time  
unsymmetry it is advisable
to choose a time direction that is to be maintained throughout the exposition. 
In principle the choice of time direction is arbitrary. 
For this paper 
we choose the time direction that is the opposite of the time direction that was chosen  
in \cite{Kr6}.  
However,
we do not change the definition of property (D). 
Instead, we introduce  a notion  of
"$\lambda$-synchronization"  that is equivalent to the time symmetric  
opposite of property (D).
Property (D), and therefore also $\lambda$-synchronization, is an invariant of topological conjugacy
(\cite[Proposition 4.3]{Kr6}).

$\lambda$-graph systems were introduced in \cite{Ma}.
We will recall their definition in Section 3. 
There is a one-to-one correspondence between separated one right  
resolving $\lambda$-graph systems and compact Shannon graphs that present a subshift 
$X\subset \Sigma^{\Bbb Z}$ (see \cite{KM}). 
For a subshift $X\subset \Sigma^{\Bbb Z}$ with property (D) 
there was  
constructed in \cite{Kr2000}  a compact Shannon graph
${\cal G}_D(X)$ that presents X, 
and that is invariantly associated to $X$,
and that generalizes the right Fischer cover \cite{Fis}. 
Bypassing the
compact Shannon graph ${\cal G}_D(X)$ we give in Section 3 
a direct construction of the
$\lambda$-graph system that corresponds to ${\cal G}_D(X)$ 
(or rather, due to the different choice of time direction, 
of its time symmetric opposite, 
which generalizes the left Fischer cover,
the subshift  $X$ now being
assumed to be $\lambda$-synchronizing ). 
We also give a direct proof  
that this $\lambda$-graph system, 
that we call the $\lambda$-synchronizing $\lambda$-graph system, 
is invariantly
associated to the $\lambda$-synchronizing subshift.

In Sections 4.1 and 4.2  we give examples of subshifts that are
$\lambda$-synchronizing and have property (D).

In Section 5 we consider the $C^*$-algebras that are obtained from the  
$\lambda$-synchronizing
$\lambda$-graph systems of $\lambda$-synchronizing subshifts.

{\it Acknowledgements.}
The insight of the referee leads to a substantial improvement of Section 4.2. 
This work was supported by Grant-in-Aid for Scientific Reserch (20540215),
Japan Societey for the Promotion of Science.

\section{$\lambda$-synchronization}

Let
$X \subset \Sigma^{\Bbb Z}$
be a subshift,
and let 
$l \in {\Bbb N}$.
We say that a word $v \in {\cal L}(X)$ 
is $l$-synchronizing if
$\Gamma^-_l(v) \subset \omega^-_l(v)$.
We denote the set of
 $l$-synchronizing words of $X$ by ${\cal S}_l(X)$.
We say that  
$X$ is $\lambda$-synchronizing 
if
for $w \in {\cal L}(X)$ and $k \in {\Bbb N}$
there is a word 
$v \in {\cal S}_k(X) \cap \Gamma^+(w)$.
\begin{lem}
For a subshift $X\subset \Sigma^{\Bbb Z}$,
the following are equivalent:
 \begin{enumerate}
\renewcommand{\labelenumi}{(\roman{enumi})}
\item $(X,S_{\Sigma}^{-1})$ has property (D).
\item $X$ is $\lambda$-synchronizing.
\item For $b\in {\cal L}(X)$ there exists an $a\in {\cal L}(X)$ such that  
$b\in \omega^-(a)$.
\end{enumerate}
\end{lem}
\begin{pf}
(i) $\Rightarrow$ (iii): Asume $(i)$ and let $b\in {\cal L}(X)$. For  
$\sigma \in \Gamma^+_1(b)$ there exists by (i) a $c \in {\cal L}(X)$  
such that one has for $a = \sigma c$ that $b\in \omega^-(a)   $.

(iii) $\Rightarrow$ (i):
Assume (iii)  and let $\sigma\in \Sigma$ and 
$a\in \Gamma^+(\sigma).  
$ 
By (iii) there exists a $b\in \Gamma^+(\sigma a)$ such that  
$\sigma a \in \omega ^-(b)$, and this implies that
$\sigma \in \omega ^-(ab).$

(iii) $\Rightarrow$ (ii):
Assume (iii) and for $w \in {\cal L}(X)$ choose a 
$b_0\in {\cal L}(X)$ 
such that 
$w\in \omega^-(b_0)$. 
Let then $k \in {\Bbb N}$, 
set
$$
Q = \card (\Gamma^-_k(b_0)),
$$
and order the set $\Gamma^-_k(b_0)$, 
writing
$$
\Gamma^-_k(b_0) = \{ c_q : 1 \leq q \leq Q  \}.
$$
Applying (ii) 
and \cite[Lemma 2.3]{Kr6}, 
one has an inductive procedure  
that yields an $R \in {\Bbb N}$, 
and indices $q_r, 1 \leq  r \leq R,$ 
such that
$$
1 \leq q_{r-1} < q_r \leq Q, \quad 1 < r \leq R,
$$
together with words  $b_r\in {\cal L}(X), 1 \leq r \leq R$,    
such that
$$
c_{q_r} \in \omega^-_k((b_s)_{0 \leq  s  \leq r})) ,
$$
and
$$
q_r = \min \{ q > q_{r-1}   : c_q \notin
\omega^-_k((b_s)_{0 \leq  s <r})) \}, \quad 1 < r \leq R,
$$
Then
$$
\{c_{q_r} : 0 \leq    r \leq R \} 
= 
\Gamma^-_k((b_r)_{0 \leq  r \leq  R})) 
= \omega^-_k((b_r)_{0 \leq  r \leq R})),
$$
and
$$
w \in \Gamma^-((b_r)_{0 \leq  r \leq R}),
$$
and (ii) is shown.

(ii) $\Rightarrow$ (iii):
Assume (ii), 
let $b \in {\cal L}(X)$, 
and let $K$ be the length of  $b$. 
By (ii) there exists an $a \in {\cal S}_K(X)\cap\Gamma^+_K(b)$  
such that $b\in \omega^-(a)$.
\end{pf}

\section{$\lambda$-synchronizing $\lambda$-graph systems}

In this section we recall the description of 
$\lambda$-graph systems and related invariants to define 
$\lambda$-synchronizing $\lambda$-graph systems.
  
Notions of 
$\lambda$-graph system and
symbolic matrix system have been introduced in \cite{Ma}.
They are presentations of subshifts
and generalizations  of 
finite labeled graphs and  symbolic matrices
 respectively.
A $\lambda$-graph system 
$ {\frak L} = (V,E,\lambda,\iota)$
over $\Sigma$
consists of a vertex set 
$V = V_0\cup V_1\cup V_2\cup\cdots$, an edge set 
$E = E_{0,1}\cup E_{1,2}\cup E_{2,3}\cup\cdots$, 
a labeling
$\lambda: E \rightarrow \Sigma$
and a surjective map
$\iota_{l,l+1}: V_{l+1} \rightarrow V_l$ for each
$l\in \Zp.$
An edge $e \in E_{l,l+1}$ has its source vertex $s(e)$ in $V_l,$
its terminal vertex $t(e)$ in $V_{l+1}$ 
and its label $\lambda(e)$ in $\Sigma$.
It is then required that there exists an edge in $E_{l,l+1}$
with label $\alpha$ and its terminal is  $v \in V_{l+1}$
 if and only if 
 there exists an edge in $E_{l-1,l}$
with label $\alpha$ and its terminal is $\iota(v) \in V_{l}.$
For 
$u \in V_{l-1}$ and
$v \in V_{l+1},$
we put
\begin{align*}
E^{\iota}(u,v)
& = \{e \in E_{l,l+1} \mid t(e) = v, \iota(s(e)) = u \},\\
E_{\iota}(u,v)
& = \{e \in E_{l-1,l} \mid s(e) = u, t(e) = \iota(v) \}.
\end{align*}
Then there exists a bijective correspondence
between 
$
E^{\iota}(u,v)
$
and
$
E_{\iota}(u,v)
$
that preserves labels
for all pairs $(u,v) \in V_{l-1}\times V_{l+1}$ of vertices.
This property is called the local property of the $\lambda$-graph system.

A symbolic matrix  system $(\M,I)$ consists of
a sequence  of pairs $(\M_{l,l+1}, I_{l,l+1}),$
$ l \in \Zp,$
of rectangular symbolic matrices
$\M_{l,l+1}$
and
rectangular $\{0,1\}$-matrices
$I_{l,l+1}$, 
where $\Zp$ denotes the set of all 
nonnegative integers.
Both the matrices $\M_{l,l+1}$ and $I_{l,l+1}$ have the same size
for each $l \in \Zp.$
The column size of $\M_{l,l+1}$ is the same as the row size of 
$
\M_{l+1,l+2}.
$
They satisfy the following commutation relations as symbolic matrices
\begin{equation}
I_{l,l+1} \M_{l+1,l+2} = \M_{l,l+1}I_{l+1,l+2}, 
\qquad
l \in \Zp. 
\end{equation}
We further assume  that
for $i$ 
there exists $j$ such that
the $(i,j)$-component $I_{l,l+1}(i,j) =1,$
and for  $j$ 
there uniquely exists $i$ such that
$I_{l,l+1}(i,j) =1$.

For a symbolic matrix system  $(\M,I)$,
the labeled edges from a vertex $v_i^l \in V_l$ 
to a vertex  $v_j^{l+1}\in V_{l+1}$ are given by the symbols appearing in the
$(i,j)$-component $\M_{l,l+1}(i,j)$
of $\M_{l,l+1}$.
The matrix $I_{l,l+1}$ 
defines a surjection $\iota_{l,l+1}$ 
from $V_{l+1}$ to $V_l$ for each $l \in \Zp.$
By this observation, 
the symbolic matrix systems and the $\lambda$-graph systems 
 are the same objects.
 We say that a $\lambda$-graph system ${\frak L}$ 
 presents a subshift $X$ if 
 the set ${\cal L}(X)$ of admissible words of $X$
coincides with the set of finite label sequences
 appearing in the labeled Bratteli diagram for 
${\frak L}$.

For a symbolic matrix system
$(\M,I)$, 
let
$M_{l,l+1}$ be the nonnegative rectangular matrix obtained from 
$\M_{l,l+1}$ by setting all the symbols equal to $1$ 
for each $l \in \Zp$.
Then the resulting pair
$(M,I)$ satisfies the  following relations by (3.1)
\begin{equation}
I_{l,l+1} M_{l+1,l+2} = M_{l,l+1}I_{l+1,l+2}, 
\qquad
l \in \Zp. 
\end{equation}
We call
$(M,I)$
the nonnegative matrix system
for $(\M,I)$.

For topological Markov shifts and sofic shifts, 
several topological conjugacy invariants and flow equivalence invariants,
such as dimension groups (\cite{Kr2}, \cite{Kr3}) and Bowen-Franks groups 
(\cite{BF}, \cite{Fr}) 
have been defined by using underlying matrices
(cf. \cite{BK}, \cite{Na}, \cite{Ki}, \cite{LM}). 
These invariants  have been generalized 
to nonnegative matrix systems in \cite{Ma}.
For a nonnegative matrix system $(M,I)$,
let $m(l)$ be the row size of the matrix $I_{l,l+1}$ for each $l\in \Zp$.
Let 
${\Bbb Z}_{I^t}$
be the abelian group defined by the inductive limit
$
{\Bbb Z}_{I^t} = 
\varinjlim \{ I_{l,l+1}^t : {\Bbb Z}^{m(l)} \rightarrow {\Bbb Z}^{m(l+1)} \}.
$
The sequence $ M_{l,l+1}^t, l \in \Zp$ of the transposes of $M_{l,l+1}$
 naturally acts on
 ${\Bbb Z}_{I^t}$ by the relation (3.2),  
 that is denoted by $\lambda_{(M,I)}$.
 The K-groups for $(M,I)$  have been  defined as:
\begin{equation*}
K_0(M,I)  = {\Bbb Z}_{I^t} / (\id - \lambda_{(M,I)}){\Bbb Z}_{I^t},\qquad
K_1(M,I)   = \Ker(\id - \lambda_{(M,I)} ) 
              \text{ in } {\Bbb Z}_{I^t}.
\end{equation*}
Set
the inductive limits
$
{\Bbb Z}_{I^t}^{+}  = \underset{l}{  
\underset{\longrightarrow}{\lim} }
 \{ I^{t}_{l,l+1}: {\Bbb Z}^{m(l)}_+ 
            \longrightarrow {\Bbb Z}^{m(l+1)}_+ \}
$
of positive cones.
We put ${\Bbb Z}_{I^t}(k)={\Bbb Z}_{I^t}, k\in {\Bbb N}$
and consider the inductive limits:
\begin{align*}
\Delta_{(M,I)} & = 
\underset{k}{
\underset{\longrightarrow}{\lim}} \{ \lambda_{(M,I)}:{\Bbb Z}_{I^t}(k) 
                    \longrightarrow {\Bbb Z}_{I^t}(k+1) \}, \\
\Delta^+_{(M,I)} & = 
\underset{k}{\underset{\longrightarrow}{\lim}}\{ \lambda_{(M,I)}:{\Bbb Z}^+_{I^t}(k) 
                    \longrightarrow {\Bbb Z}^+_{I^t}(k+1) \}.
\end{align*}
The ordered group
$(\Delta_{(M,I)},\Delta^+_{(M,I)})$
is called the dimension group for $(M,I)$.
The map
$
\delta_{(M,I)} :{\Bbb Z}_{I^t}(k) 
                    \rightarrow {\Bbb Z}_{I^t}(k+1) 
$
defined by
$\delta_{(M,I)}([X,k]) = ([X,k+1])$
yields an automorphism
on 
$(\Delta_{(M,I)},\Delta^+_{(M,I)}).$ 
The triple
$(\Delta_{(M,I)},\Delta^+_{(M,I)}, \delta_{(M,I)})$
is named the dimension triple for $(M,I).$ 
We set the projective limit of the abelian group as 
$
{\Bbb Z}_I = \underset{\longleftarrow}{\lim} \{I_{l,l+1}: {\Bbb Z}^{m(l+1)} \longrightarrow 
                                      {\Bbb Z}^{m(l)} \}.
$ 
The sequence $M_{l,l+1}, l \in \Zp$ 
 acts on ${\Bbb Z}_I$
as an endomorphism that we denote by
$M$.
The identity on ${\Bbb Z}_{I}$ 
is denoted by $I$.
The Bowen-Franks groups have been formulated as:  
\begin{equation*}
BF^0(M,I) = {\Bbb Z}_I /(I-M){\Bbb Z}_I, 
\qquad
BF^1(M,I) = \Ker(I-M) \text{ in }{\Bbb Z}_I.
\end{equation*}
Both the pairs $K_*(M,I), BF^*(M,I)$ for $* =0,1$ are invariant under shift 
equivalence in nonnegative matrix systems (\cite{Ma}).

In \cite{Ma2002a}, 
the $C^*$-algebra ${\cal O}_{\frak L}$
associated with a $\lambda$-graph system ${\frak L}$ has been introduced.
These $C^*$-algebras are generalizations of the Cuntz-Krieger algebras 
and the $C^*$-algebras associated with subshifts.
They are universal unique concrete $C^*$-algebras generated 
by finite families of partial isometries 
and sequences of projections subject to certain operator relations  
encoded by structure of the $\lambda$-graph systems.
The $C^*$-algebra 
${\cal O}_{\frak L}$ has a natural one-parameter group action  
called gauge action.
Its fixed point algebra denoted by 
${\cal F}_{\frak L}$
becomes an AF-algebra.
Let $(M,I)$ be the nonnegative matrix system 
for the symbolic matrix system of the $\lambda$-graph system
${\frak L}$.
The following relations hold:
\begin{align}
K_0({\cal O}_{\frak L}) & = K_0(M,I), \qquad
K_1({\cal O}_{\frak L}) = K_1(M,I),\\
\Ext^1({\cal O}_{\frak L}) & = BF^0(M,I),\qquad
\Ext^0({\cal O}_{\frak L}) = BF^1(M,I), \\
(K_0({\cal F}_{\frak L}),& K_0({\cal F}_{\frak L})^+)  =
(\Delta_{(M,I)}, \Delta^+_{(M,I)})
\end{align}
where 
$\Ext^1({\cal O}_{\frak L})= \Ext({\cal O}_{\frak L})$
and
$\Ext^0({\cal O}_{\frak L})= \Ext({\cal O}_{\frak L}\otimes C_0({\Bbb R}).$

All the groups above are invariant under shift equivalence of symbolic matrix systems,
so that they yield topological conjugacy unvariants of subshifts by
taking canonical $\lambda$-graph systems.

The $\lambda$-entropy  
$h_{\lambda}(\frak L)$
of a $\lambda$-graph system ${\frak L}$
was introduced in \cite{KM3}.
The $\lambda$-entropy
 measures the growth rate of the cardinalities 
${}^{\sharp}V_l, l \in \Zp$
of the vertex sets 
$\{ V_l \}_{l \in \Zp}$.
The volume entropy 
$
h_{vol}({\frak L})
$
of a $\lambda$-graph system ${\frak L}$
was introduced in \cite{Ma2005c}.
Denote by $P_l({\frak L})$ the set of all labeled paths starting at
a vertex in  
$V_0$ and terminating at a vertex in $V_l$.
The volume entropy 
measures the growth rate of the cardinalities ${}^{\sharp}P_l(\frak L)$ 
of the labeled paths $P_l({\frak L})$.
Both the entropic quantities are
 invariant under shift equivalence of $\lambda$-graph systems,
so that they yield a topological conjugacy invariants of subshifts.

\medskip

For $\mu, \nu \in {\cal L}(X)$,
if $\Gamma_l^-(\mu) = \Gamma_l^-(\nu)$,
we say that $\mu$ is $l$-past equivalent to $\nu$
and write it as  
$\mu\underset{l}{\sim}\nu$.
\begin{lem}
Let $X$ be a $\lambda$-synchronizing subshift.
Then we have
\begin{enumerate}
\renewcommand{\labelenumi}{(\roman{enumi})}
\item 
For $l \in {\Bbb N}$ and $\eta \in {\cal L}_n(X)$,
there exists
$\mu \in {\cal S}_l(X)$
such that $\eta \in \Gamma_l^-(\mu)$. 
\item
For $\mu \in {\cal S}_l(X)$,
there exists
$\mu' \in {\cal S}_{l+1}(X)$
such that 
$\mu\underset{l}{\sim}\mu'$.
\item
For $\mu \in {\cal S}_l(X)$,
there exist
$\beta \in \Sigma$
and
$\nu \in {\cal S}_{l+1}(X)$
such that 
$\mu\underset{l}{\sim}\beta\nu$.
\end{enumerate}
\end{lem}
\begin{pf}
(i)
 This follows from  $\lambda$-synchronization.

(ii)
For $\mu \in {\cal S}_l(X)$ with
$|\mu | = K$,
put 
$k = K + l+1 > K$.
As $X$ is $\lambda$-synchronizing,
there exists $\nu \in {\cal S}_k(X)$
such that
$\mu \nu \in {\cal S}_{k - K}(X)$.
Put $\mu'=\mu \nu  \in {\cal S}_{l+1}(X)$.
As $\mu \in {\cal S}_l(X)$, one sees that 
$\Gamma_l^-(\mu) = \Gamma_l^-(\mu \nu)$
so that 
$\mu\underset{l}{\sim}\mu'$.

(iii)
For $\mu \in {\cal S}_l(X)$ with
$\mu = \mu_1 \cdots \mu_K$,
put
$ k = K +l >K$.
As $X$ is $\lambda$-synchronizing,
there exists $\omega \in {\cal S}_k(X)$
such that
$\mu \omega \in {\cal S}_{k - K}(X)$.
Set
$\beta = \mu_1$
and
$\nu = \mu_2 \cdots \mu_K \omega$.
Since $\omega \in {\cal S}_k(X)$,
one has 
$\nu \in {\cal S}_{k-(K-1)}(X)$
so that 
$\nu \in {\cal S}_{l+1}(X)$.
As 
$\Gamma_l^-(\mu) = \Gamma_l^-(\mu \omega)$,
one sees
that 
$\mu\underset{l}{\sim}\beta \nu$.
\end{pf}

For a $\lambda$-synchronizing subshift $X$ over $\Sigma$,
we will introduce the $\lambda$-synchronizing $\lambda$-graph system
\begin{equation*}
{\frak L}^{\lambda(X)} 
=(V^{\lambda(X)}, E^{\lambda(X)}, \lambda^{\lambda(X)}, \iota^{\lambda(X)})
\end{equation*}
in the following way.
Let
$V_l^{\lambda(X)}$ be the $l$-past equivalence classes of
${\cal S}_l(X)$.
We denote by $[\mu]_l$ 
the  $l$-past equivalence class of $\mu \in {\cal S}_l(X)$.
For
$\nu \in {\cal S}_{l+1}(X)$ and $\alpha \in \Gamma_1^-(\nu)$,
define an edge with label $\alpha$
 from
$[\alpha \nu]_l\in V_l^{\lambda(X)}$ to
$[\nu]_l \in V_{l+1}^{\lambda(X)}$.
We denote the set of these edges by $E^{\lambda(X)}_{l,l+1}$.
Since
${\cal S}_{l+1}(X) \subset {\cal S}_l(X)$,
we have a natural map
$
[\mu]_{l+1} \in V_{l+1}^{\lambda(X)} \longrightarrow 
[\mu]_l \in V_l^{\lambda(X)} 
$
that we denote by
$\iota^{\lambda(X)}_{l,l+1}$.
\begin{prop}
${\frak L}^{\lambda(X)} 
=(V^{\lambda(X)}, E^{\lambda(X)}, \lambda^{\lambda(X)}, \iota^{\lambda(X)})
$
is a $\lambda$-graph system that presents $X$.
\end{prop}
\begin{pf}
We will first  show the local property of $\lambda$-graph systems
(cf. \cite{Ma2002a}, \cite{Ma2005b}).
For $[\mu]_l \in V_l^{\lambda(X)}$
and
$[\nu]_{l+2} \in V_{l+2}^{\lambda(X)}$
with
$\mu \in {\cal S}_l(X),\nu \in {\cal S}_{l+2}(X)$,
suppose that there exists a labeled edge 
from
$[\mu]_l$ to $[\nu]_{l+1}$
labeled $\alpha \in \Sigma$.
It follows that
$\alpha \nu \underset{l}{\sim}\mu$.
Hence there exists an edge from
$[\alpha \nu]_{l+1}$ 
to
$[\nu]_{l+2}$
labeled $\alpha$
and a $\iota$-map from
$[\alpha \nu]_{l+1}$ to 
$[\alpha \nu]_l$.

On the other hand,
suppose that there exists an
$\iota$-map from
$[\omega]_{l+1}$
to
$[\mu]_l$
and an edge from
$[\omega]_{l+1}$
to
$[\nu]_{l+2}$
labeled $\alpha$.
Hence
$\omega\underset{l+1}{\sim}\alpha\nu$.
Since
$\iota^{\lambda(X)}([\alpha\nu]_{l+1})= 
[\alpha\nu]_l$,
one sees that
$[\mu]_l =
[\alpha\nu]_l$
so that
$
\mu\underset{l}{\sim}
\alpha\nu$.
Hence there exists an edge from
$[\mu]_l$
to
$[\nu]_{l+1}$
labeled $\alpha$.
Therefore the local property of $\lambda$-graph systems holds.

By $\lambda$-synchronization,
for any admissible word $\eta \in {\cal L}(X)$ and $k \ge l$,
there exists $\nu \in {\cal S}_k(X)$ such that
$\eta \nu \in {\cal S}_{k-l}(X)$.
This implies that
there exists a path labeled $\eta$ in ${\frak L}^{\lambda(X)}$
from the vertex 
$[\eta \nu]_{k-l} \in V_{k-l}^{\lambda(X)}$
to the vertex
$[\nu]_l \in V_l^{\lambda(X)}$
so that 
${\frak L}^{\lambda(X)}$
is a $\lambda$-graph system that presents $X$.
\end{pf}

As in \cite{Ma},
there is a canonical construction of a $\lambda$-graph system
for an arbitrary subshift $X$.
The constructed $\lambda$-graph system is called 
the canonical $\lambda$-graph system for $X$
and denoted by ${\frak L}^X$.
For $l \in \Zp$, 
the vertex set $V_l^X$ of ${\frak L}^X$
is defined by the $l$-past equivalence classes 
$\Gamma_l^-(x)$ 
of the right infinite sequences $x \in X^+$.
For a symbol $\alpha \in \Sigma$,
an edge labeled $\alpha$ from 
$\Gamma_l^-(\alpha x)$ 
to
$\Gamma_{l+1}^-( x)$
is defined if
$\alpha x \in X^+$.
The natural inclusions
$\Gamma_{l+1}^-( x) \subset \Gamma_l^-( x)$
give rise to the $\iota$-map.
\begin{cor}
$
{\frak L}^{\lambda(X)} 
$
is a predecessor-separated, left-resolving $\lambda$-graph subsystem
of the canonical $\lambda$-graph system of $X$.
\end{cor}
\begin{pf}
One checks  that 
$
{\frak L}^{\lambda(X)} 
$
is  predecessor-separated and left-resolving.
Let
${\frak L}^{X} 
=(V^X, E^X,\lambda^X,\iota^X)
$
be the canonical $\lambda$-graph system
for
$X$.
 For $\mu \in {\cal S}_l(X)$
and $x \in X^+$ with
$\mu \in \Gamma_l^-(x)$,
one sees that 
$\Gamma_l^-(\mu) = 
\Gamma_l^-(\mu x).
$ 
Hence
$[\mu]_l$ can be regarded as a vertex of $V_l^X$
so that 
the vertex set
$V^{\lambda(X)}_l$ can be regarded as a subset of $V_l^X$.
Similarly
the edge set
$E_{l,l+1}^{\lambda(X)}$ 
can be regarded as a subset of $E_{l,l+1}^X$.
The $\iota$-map $\iota^{\lambda(X)}$ of 
${\frak L}^{\lambda(X)}$
is compatible to
that of ${\frak L}^X$.
It follows  that   
$
{\frak L}^{\lambda(X)} 
$
is a $\lambda$-graph subsystem
of ${\frak L}^X$.
\end{pf}

For a synchronizing subshift $X$,
the canonical synchronizing $\lambda$-graph system
${\frak L}^{S(X)}$ of $X$
has been introduced in 
\cite{KM}.
\begin{prop}
Let $X$ be a synchronizing  subshift
and
${\frak L}^{S(X)}$ the canonical synchronizing $\lambda$-graph system of $X$.
Then 
${\frak L}^{S(X)}$ is isomorphic to 
the 
$\lambda$-synchronizing $\lambda$-graph system 
${\frak L}^{\lambda(X)}$
of $X$.
\end{prop}

For a $\lambda$-synchronizing subshift 
$X$,
denote by
$(\M^{\lambda(X)}, I^{\lambda(X)})$
the symbolic matrix system for
the $\lambda$-graph system
$
{\frak L}^{\lambda(X)} 
$ (see \cite{Ma}).
\begin{prop}
Let 
$X, X'$ be $\lambda$-synchronizing subshifts.
If $X$
is topologically conjugate to $X'$,
then 
$(\M^{\lambda(X)}, I^{\lambda(X)})$
is strong shift equivalent to
$(\M^{\lambda(X')}, I^{\lambda(X')})$.
\end{prop}
\begin{pf}
By \cite{Na}, we may assume that 
the subshifts 
$X, X'$ are bipartitely related to each other.
This means that there exists alphabets $C, D $
and a bipartite subshift 
$\widehat{X}$ over $C\cup D$ such that
there exist
subshifts 
$X_1$ over $CD$, $X_2$ over $DC$
and
one block conjugacies
$\varphi_1 : X \longrightarrow X_1$,
$\varphi_2 : X' \longrightarrow X_2$
such that
\begin{equation*}
{\widehat{X}}^{[2]} = X_1 \cup X_2,
\end{equation*}
where
${\widehat{X}}^{[2]}$ is the 2-block shift of $X$.
We will show that 
the $\lambda$-graph systems
$
{\frak L}^{\lambda(X)} 
$
and
$
{\frak L}^{\lambda(X')} 
$
are bipartitely related to each other.
We note that
\begin{equation}
{\cal S}_{2k}({\widehat{X}}) = \varphi_1({\cal S}_k(X)) \cup \varphi_2({\cal S}_k(X')). 
\label{eqn:bip}
\end{equation}
For
$\eta = (\eta_1,\dots,\eta_l) \in {\cal L}_l(\widehat{X})$,
assume that  $\eta_1 \in C$.

Case 1: $\eta_l \in C$.

Hence $l $ is odd, $l = 2m-1$ for some 
$m \ge 1$.
Take $d \in D$ such that 
$
\eta_1\cdots\eta_l d \in {\cal L}_{l+1}(\widehat{X}) = {\cal L}_{2m}(\widehat{X})
$
so that
$\varphi_1^{-1}(\eta_1\cdots\eta_l d) \in {\cal L}_{m}(X).$
For
$k \ge l$,
take 
$k'$ such as 
$k= 2 k' -1$ if $k$ is odd,
and
$k= 2 k'$ if $k$ is even.
Hence $k' +1 \ge m$.
Since $X$ is $\lambda$-synchronizing,
there exists $\nu \in {\cal L}_{k' +1}(X)$
such that
$\varphi_1^{-1}(\eta_1\cdots\eta_l d)\nu \in {\cal S}_{k'+1-m}(X).$
As $2(k'+1-m)=2k'-2(m-1) \ge k-l$,
one has
$\eta_1\cdots\eta_l d \varphi_1(\nu) \in {\cal S}_{k-l}(\widehat{X}).$

Case 2: $\eta_l \in D$.

Hence $l $ is even, $l = 2m$ for some 
$m \ge 1$.
One sees
$\varphi_1^{-1}(\eta_1\cdots\eta_l ) \in {\cal L}_{m}(X).$
For
$k \ge l$,
take 
$k'$ such as 
$k= 2 k' -1$ if $k$ is odd,
and
$k= 2 k'$ if $k$ is even.
Hence $k' \ge m$.
Since $X$ is $\lambda$-synchronizing,
there exists $\omega \in {\cal L}_{k'}(X)$
such that
$\varphi_1^{-1}(\eta_1\cdots\eta_l d)\omega\in {\cal S}_{k'-m}(X).$
As $2k'- 2m \ge k-l$,
one has
$\eta_1\cdots\eta_l d \varphi_1(\omega)\in  {\cal S}_{k-l}(\widehat{X}).$

Thus the bipartite subshift
$\widehat{X}$ is $\lambda$-synchronizing.
The equality \eqref{eqn:bip}
implies 
\begin{equation*}
V_{2k}^{\lambda(\widehat{X})} = 
V_{k}^{\lambda(X)} \cup V_{k}^{\lambda(X')}.
\end{equation*}
One then easily sees that 
the $\lambda$-graph systems
${\frak L}^{\lambda(X)}$
and
${\frak L}^{\lambda(X')}$
are bipartite pair in the sense of \cite{Na}.
Hence
$(\M^{\lambda(X)}, I^{\lambda(X)})$
is strong shift equivalent to
$(\M^{\lambda(X')}, I^{\lambda(X')})$
by \cite{Ma}.
\end{pf}

For the symbolic matrix system
$(\M^{\lambda(X)}, I^{\lambda(X)})$,
denote by
$(M^{\lambda(X)}, I^{\lambda(X)})$,
its nonnegative matrix system.
Set
\begin{align}
K_i^\lambda(X) & = K_i(M^{\lambda(X)},I^{\lambda(X)}),\qquad i=0,1 \\ 
BF^i_\lambda(X) & = BF^i(M^{\lambda(X)},I^{\lambda(X)}),\qquad i=0,1 \\ 
(\Delta^\lambda(X), & \Delta_+^\lambda(X)) 
= (\Delta_{(M^{\lambda(X)},I^{\lambda(X)})},\Delta^+_{(M^{\lambda(X)},I^{\lambda(X)})}),\\ 
h_\lambda(X) & = h_\lambda({\frak L}^{\lambda(X)}), 
\\
h_{vol}^\lambda(X) & = h_{vol}({\frak L}^{\lambda(X)}).
\end{align}
Since the above invariants are all shift equivalence invariants
and
since strong shift equivalence implies shift equivalence,
we have
\begin{cor}
The abelian groups 
$
K_i^\lambda(X), BF^i_\lambda(X), i=0,1,
$
the ordered abelian group
$
\Delta^\lambda(X)
$
and the entropic quantities 
$h_\lambda(X),
h_{vol}^\lambda(X)
$
 are all invariant under topological conjugacy of $\lambda$-synchronizing subshifts.
\end{cor}

\section{Examples}

{\bf 1. Dyck shifts and Motzkin shifts.}

Starting from the Dyck shifts and the Motzkin shifts, we describe, in  
increasing generality, classes of subshifts that have been observed to  
have property (D).
 These subshifts are also $\lambda$-synchronizing,  
since these classes are closed under taking inverses. The descriptions  
involve finite directed graphs. The mapping that assigns to a path in  
a directed graph its source vertex we denote by $s$ and the mapping  
that assigns to a path in the graph its target vertex we denote by $t$.

First we recall the construction of the Dyck shifts and of the Motzkin shifts.
We denote the generators of  the Dyck inverse monoid (the polycyclic  
inverse monoid) by $e_n^- ,e_n^+,1 \le n \le N, N >1$. These  
generators satisfy the relations
\begin{equation*}
e_n^- e_n^+ =  {\bold 1}, \quad 1 \le n\le N, \qquad
e_l^- e_m^+ = 0, \quad 1 \le l, m \le N, \ l \ne m.
\end{equation*}
The Dyck shift $D_{N}$ is the subshift with alphabet
  $\{e_n^- ,e_n^+   : 1 \le n\le N\}$
and admissible words
$(e_{i})_{1 \leq i \leq I}, I \in \Bbb N,$ given by the condition
$
\prod_{1 \leq i \leq I} e_{i} \neq  0.
$
The Motzkin shift  $M_{N }$ is the subshift with alphabet
$\{e_n^- ,e_n^+   : 1 \le n\le N\}\cup \{ \bold 1  \} $
and admissible words
$(e_{i})_{1 \leq i \leq I}, I \in \Bbb N,$  also given by the condition
$
\prod_{1 \leq i \leq I} e_{i} \neq  0.
$

The Dyck shifts belong to the class of Markov-Dyck shifts and the  
Motzkin shifts belong to the class of Markov-Motzkin shifts. 
To recall the construction of the Markov-Dyck shifts and of the Markov-Motzkin shifts, 
let there be given an irreducible finite directed graph
with vertex set ${\mathcal V}$ and edge set
${\mathcal E}$. Let $({\mathcal V}, {\mathcal E}^-  )$ 
be
a copy of  $({\mathcal V}, {\mathcal E}  )$.
Reverse the directions of the edges in ${\cal E}$
to obtain the reversed graph
of the graph  $({\mathcal V}, {\mathcal E} )$
with vertex set ${\cal V}$ and edge set
${\cal E}^+$.
With idempotents $P_v, v \in {\cal V},$
the set
${\cal E}^- \cup \{P_v:v \in {\cal V}\}\cup {\cal E}^+$
is the generating set of the graph inverse semigroup ${\mathcal  
S}_{{\mathcal V}, {\mathcal E}}$ of the directed graph 
$({\mathcal V},  {\mathcal E})$, 
where, besides
$P_u^2 = P_u, v \in {\cal V},$ 
the relations are
$$
P_uP_w = 0, \qquad   u, w \in {\cal V}, \ u \ne w,
$$
\begin{equation*}
f^-g^+
=
\begin{cases}
P_{s(f)}, & (f = g),\\
0,  &  (f \ne g, \ f, g \in {\cal E} ),
\end{cases}
\end{equation*}
and
\begin{equation*}
P_{s(f)} f^- = f^- P_{r(f)}, \qquad 
P_{r(f)} f^+ = f^+ P_{s(f)}, \qquad f \in {\cal E}.
\end{equation*}
The Markov-Dyck shift  of the graph $({\mathcal V}, {\mathcal E})$  
is the subshift with alphabet
${\cal E}^-\cup {\cal E}^+$
with  admissible words $(e_{i})_{1 \leq i \leq I}, I \in {\Bbb N},$
given  
by the condition
$
\prod_{1 \leq i \leq I} e_{i} \neq  0.
$
The Dyck shift $D_N$ arises in this way from the single vertex graph  
with $N$ loops at its vertex.
The Markov-Motzkin shift  of the graph $({\mathcal V}, {\mathcal E}   
)$ is the subshift with alphabet
${\cal E}^-\cup \{  P_v: v \in {\mathcal V} \} \cup {\cal E}^+$
with  admissible words
$(e_{i})_{1 \leq i \leq I}, I \in \Bbb N,$ also  given by the condition
$
\prod_{1 \leq i \leq I} e_{i} \neq  0.
$ 
The Motzkin shift $M_N$ arises in this way from the single vertex  
graph with $N$ loops at its vertex.

Following \cite{HI}, in \cite{HIK} a necessary and sufficient condition was given  
for the existence of an embedding of
an irreducible subshift of finite type into target subshifts that were  
taken from a class of  
$\lambda$-synchronizing subhifts with property (D). 
This class contains the Markov-Dyck and Markov-Motzkin shifts.  
To recall the construction of this class, let there be given,
besides the finite irreducible directed graph $({\mathcal V}, {\mathcal E})$,  
another finite irreducible directed graph with vertex set $\Omega$ 
and edge set $\Sigma$. 
Denote by 
${\mathcal  S}^-_{\mathcal V,  \mathcal E}$
(resp.  $ {\mathcal  S}^+_{\mathcal V,  \mathcal E}$) 
the semigroup that  
is generated by 
$\{e^-: e \in
{\mathcal E  }\}$
(resp. $\{e^-: e \in {\mathcal E  }\}$), 
and let $\lambda$ be a labeling map that assigns to every edge 
$\sigma \in\Sigma $ 
a label
$\lambda(\sigma) 
\in
     {\mathcal  S}^-_{ \mathcal V,  \mathcal E} 
\cup \{ P_v :v \in {\mathcal V}\} 
\cup {\mathcal  S}^+_{ \mathcal V,  \mathcal E} $,
  and extend the mapping $\lambda$  to all finite paths 
$ ( \sigma_i)_{1 \leq i \leq I}$ 
in the graph $(\Omega,\Sigma)$ 
by
$$
\lambda (( \sigma_i )_{1 \leq i \leq I}) = \prod_{1 \leq i \leq I}\lambda(\sigma_i ).
$$
For $v\in {\mathcal V}$ 
let $\Omega_v$ denote the set of $\omega \in \Omega$ such  
that there exists a path $a$ in the graph $(\Omega,\Sigma)$ such that  
$s(a) = t(a) = \omega$ and $\lambda (\omega) = P_v$. 
Assume that
$\Omega_v \not= \emptyset, v \in {\mathcal V},$ 
and that
$\{\Omega_v : v \in {\mathcal V}  \} $ 
is a partition of $ \Omega$. 
Also assume that for every edge that enters an
$\omega  \in  \Omega_v$ one has $\lambda (\omega)P_v \not=0 $ and  
for every edge $\sigma$ that leaves  an  
$\omega  \in  \Omega_v$ 
one has 
$P_v \lambda (\omega)\not= 0, v \in {\mathcal V}$. 
Assume that for  
$u, w \in {\mathcal V}$ 
and $g\in {\mathcal S}_{{\mathcal V}, {\mathcal E}}$ 
such that 
$ P_ug P_w \not= 0$ 
there exists a path $a$ in the graph such that  
$s(a) = u, t(a) =w $ and $\lambda (a) = g$.
We define a subshift $X(\mathcal V,  \mathcal E, \lambda)$ 
as the subsystem of the edge
shift of $( \mathcal V,  \mathcal E)$ with allowed words the finite  
paths $b$ in $G$ such that
$\lambda(b) \not = 0$.  The class of subshifts of the form 
$X({\mathcal V},  \mathcal E, \lambda)$ is closed under taking inverses and one  
checks that the subshifts $X(\mathcal V,  \mathcal E, \lambda)$ have  
Property (D). 
In fact, they have stronger synchronization properties  
as described in \cite{Kr8}. 
The $\lambda$-synchronizing $\lambda$-graph  
system of a subshift $X(\mathcal V,  \mathcal E, \lambda)$ is given by  
its Cantor horizon as described for the Dyck shift in \cite{KM2} 
and for the  
Motzkin shift in \cite{Ma4}.
While maintaining sychronization properties one has a more general  
construction that goes beyond the class of subshifts of the form  
$X(\mathcal V,  \mathcal E, \lambda)$, where the graph inverse  
semigroup is replaced by a more general type of semigroup \cite{Kr5}.


\medskip

{\bf 2. Substituion dynamical systems}.

For the theory of substitution dynamical systems see \cite{Que}.

In proving  the next theorem we follow \cite[Example 3.6]{CE}.

\begin{thm}
The  substitution dynamical system of a primitive substitution is
$\lambda$-synchronizing.
\end{thm}
\begin{pf}
For a substitution minimal system of a primitive  substitution  
$X\subset \Sigma^{\Bbb Z}$ 
and for $k \in {\Bbb N}$, 
let $A^+(k)$ denote  
the set of $x^+ \in X_{[0, \infty )}$ 
such that 
$ \card(\Gamma ^-_k(x^+)) > 1$. 
$A^+(1)$ is a finite set \cite[Section 5.1.1]{Que} 
and  
therefore the sets 
$A^+(k),k \in {\Bbb N}$ 
are also finite. 
We can for  
$b \in {\cal L}_k(X)$
 choose a 
$y^+  \in X_{[0, \infty )} \setminus A^+(k)$ 
such that 
$b \in \Gamma^-_k( y^+ )$. 
There is an $n \in {\Bbb N}$ 
such that for $a= y^+_{[0, n]},$
$\card(\Gamma^-_k(a))  = 1$ and therefore 
$\{b\} = \omega^-_k (a)$.
Otherwise $\Gamma^-_k(y^+)$, which is the limit of the decreasing  
sequence $\Gamma^-_k(y^+_{[0, n]}), n \in {\Bbb N},$ 
would contain  
more than one word, contradicting $y^+ \notin A^+(k)$.
We have proved  that $X$ satisfies condition (iii) of Lemma 2.1.
\end{pf}

Note that we have also proved that the invariant probability measure  
of the substitution dynamical system of a  primitive substitution  is  
a $g$-measure in the sense of \cite{KW}. 
(We have used the time  direction that is opposite to the one in 
\cite{MW}, or in \cite{Kr7} or \cite{KW}.  
However, the situation is time symmetric.)


\medskip

{\bf 3. A coded system}.

With the alphabet 
$\Sigma = \{ 0,1,\alpha,\beta,\gamma\}$
we consider a subshift $X \subset \Sigma^{\Bbb Z}$
that has property (D) such that its inverse does not 
have property (D).
We obtain $X$ as the closure of the union of an increasing sequence 
$Y_n, n \in {\Bbb N}$
of irreducible subshifts of finite type.
This implies that $X$ is a coded system
\cite[Theorem 1]{Kr2000}.
$Y_n, n \in {\Bbb N}$ is defined by
excluding from $\Sigma^{\Bbb Z}$ the words
$$
0^m, \ m>n,
$$
and the words
$$
\beta \alpha   c \gamma 0^k \gamma, \qquad c \in \Sigma^k, \ 1 \le k \le n,
$$
as well as the words
$$
\beta \beta,\quad \beta \gamma,\quad \beta0, \quad \beta 1.
$$

We prove that 
$X$ has property (D).
For this let $a \sigma \in {\cal L}(X)$.
If $\sigma \not= \gamma $, 
then,
with $K$ the length of $a, \sigma \in \omega^+(1^K a)$.
If $\sigma =\gamma $, 
then $\sigma \in \omega^+(\alpha a)$.

Let 
$a\in \Gamma^+(\beta\alpha)$, 
and let $K$ be the length of $a$.  
Then $\gamma0^K\gamma \in \Gamma^+(\alpha a)$,
 but 
$\beta\alpha a  \gamma 0^K \gamma$ 
is not admissible for $X$, and therefore $X$ does  
not satisfy condition (iii) of Lemma 2.1.


\section{$C^*$-algebras}

%


Generalizing Condition (I) of \cite{CK},
in  \cite{Ma2005b}  $\lambda$-condition (I)
was introduced,
 which says that 
for a vertex $v$ in the $\lambda$-graph system,
there exist two distinct paths $\pi_1, \pi_2$ starting at 
$v$ 
such that they  have the same terminal vertex but different labels.

For a $\lambda$-synchronizing subshift $X$,
we say that $X$ satisfies {\it synchronizing condition} (I)
if for $l \in {\Bbb N}$ and
$\mu \in {\cal S}_l(X)$,
there exist 
$\gamma_1, \gamma_2 \in {\cal L}_K(X)$
for some $K$
and $\nu \in {\cal S}_{l+K}(X)$
such that
\begin{equation}
\gamma_1 \not= \gamma_2, 
\qquad 
\gamma_1,\gamma_2 \in \Gamma_K^-(\nu),
\qquad 
[\gamma_1 \nu]_l= [\gamma_2 \nu]_l= [\mu]_l. \label{eqn:syn(I)}
\end{equation}
We have the following lemma.
\begin{lem}
Let $X$ be a $\lambda$-synchronizing subshift.
Then 
the following conditions are equivalent:
\begin{enumerate}
\renewcommand{\labelenumi}{(\roman{enumi})}
\item $X$ satisfies synchronizing condition (I). 
\item The $\lambda$-synchronizing $\lambda$-graph system 
${\frak L}^{\lambda(X)} $ satisfies 
$\lambda$-condition (I). 
\end{enumerate}
\end{lem}
\begin{pf}
(i)$ \Rightarrow $(ii):
Suppose that 
$X$ satisfies synchronizing condition (I). 
For a vertex $v \in V_l^{\lambda(X)}$ in 
${\frak L}^{\lambda(X)}$,
take a $l$-synchronizing word $\mu \in {\cal S}_l(X)$
such that
$v = [\mu ]_l$.
By the synchronizing condition (I) of $X$,
there exist 
$\gamma_1, \gamma_2 \in {\cal L}_K(X)$
for some $K$
and $\nu \in {\cal S}_{l+K}(X)$
satisfying \eqref{eqn:syn(I)}.
This implies that there exist two paths beginning with $v$ and ending in
the vertex $[\nu]_{l+K} \in V^{\lambda(X)}_{l+K}$
whose labels are $\gamma_1, \gamma_2$.
Hence  
${\frak L}^{\lambda(X)} $ satisfies $\lambda$-condition (I).

(ii)$ \Rightarrow $(i):
Suppose that 
${\frak L}^{\lambda(X)} $ satisfies 
$\lambda$-condition (I). 
Let $\mu \in {\cal S}_l(X)$ be a $l$-synchronizing word of $X$.
By the $\lambda$-condition (I),
for the vertex $[\mu]_l \in V^{\lambda(X)}_l$,
there exist two distinct paths $\pi_1, \pi_2$ 
in ${\frak L}^{\lambda(X)}$
starting at 
$[\mu ]_l$ 
such that 
they have the same terminal vertex but different labels.
We denote by $u$ the terminal vertex.
As $u$ belongs to $V^{\lambda(X)}_{l+K}$
for some $K \in {\Bbb N}$,
one may find a $l+K$-synchronizing word 
$\nu \in {\cal S}_{l+K}(X)$ such that 
$[\nu]_{l+K} = u$.
 Denote by $\gamma_1,\gamma_2$
the labels of $\pi_1,\pi_2$
respectively.
Since $\pi_1,\pi_2$ 
begin with $[\mu]_l$ 
and end in $u$,
one sees that
$\gamma_1,\gamma_2 \in \Gamma_K^-(\nu)$
and
$
[\gamma_1 \nu]_l= [\gamma_2 \nu]_l= [\mu]_l.
$
As $\gamma_1\not= \gamma_2$,
one sees that 
$X$ satisfies synchronizing condition (I).
\end{pf}

If  $\lambda$-graph systems 
${\frak L}$ and $ {\frak L}'$ are bipartitely related 
by a bipartite
$\lambda$-graph system $\widehat{\frak L}$,
and if 
one of the $\lambda$-graph systems ${\frak L}, {\frak L}', \widehat{\frak L}$
satisfies $\lambda$-condition (I),
then so do the other two.
Hence $\lambda$-condition (I) is invariant under strong shift equivalence 
of the symbolic matrix systems that correspond to the $\lambda$-graph systems. 
Therefore by the preceding lemma,
one knows that synchronizing condition (I) is an invariant condition  under topological conjugacy of 
$\lambda$-synchronizing subshifts.

As a condition under which the $C^*$-algebra 
${\cal O}_{\frak L}$ is simple and purely infinite,
$\lambda$-irreducibility for $\lambda$-graph system ${\frak L}$
has been introduced in \cite{Ma2005b}.
A $\lambda$-graph system
${\frak L}$ is said to be $\lambda$-irreducible
if
for an ordered pair of vertices $v_i^l, v_j^l \in V_l,$ 
there exists a number $L_l(i,j) \in {\Bbb N}$ such that 
for a vertex $v_h^{l+L_l(i,j)} \in V_{l+L_l(i,j)}$ 
with
$\iota^{L_l(i,j)}(v_h^{l+L_l(i,j)}) = v_i^l,$
there exists a path $\gamma$ in ${\frak L}$ 
such that 
$
s(\gamma) = v_j^l, \ 
t(\gamma) = v_h^{l+L_l(i,j)},
$
where $\iota^{L_l(i,j)}$ means the $L_l(i,j)$-times compositions of $\iota$, 
and $s(\gamma), t(\gamma)$ denote 
the source vertex, 
the terminal vertex of $\gamma$ respectively. 
A $\lambda$-synchronizing subshift $X$ is said to be 
{\it synchronized irreducible} if
for $\mu, \nu \in {\cal S}_l(X)$, 
there exists $k_{\mu,\nu} \in {\Bbb N}$ such that
for $\eta \in {\cal S}_{l+ K_{\mu,\nu}}(X)$ 
with
$\nu \underset{l}{\sim} \eta$,
there exists $\xi \in {\cal L}_{k_{\mu,\nu}}(X)$ 
such that 
$\xi\eta \underset{l}{\sim} \mu$.
It is direct to see that 
a $\lambda$-synchronizing subshift $X$ is
synchronized irreducible
if and only if  
${\frak L}^{\lambda(X)}$ 
is $\lambda$-irreducible.

\begin{thm}
Let $X$ be a  $\lambda$-synchronizing subshift
satisfying synchronizing condition (I).
Suppose that $X$ is synchronized irreducible.
Then the  $C^*$-algebra
${\cal O}_{{\frak L}^{\lambda(X)}}$
associated with the $\lambda$-synchronizing $\lambda$-graph system 
${\frak L}^{\lambda(X)}$ is a simple $C^*$-algebra such that
\begin{align}
K_i({\cal O}_{{\frak L}^{\lambda(X)}})
& = K_i^\lambda(X) \qquad i=0,1 \\ 
\Ext^i({\cal O}_{{\frak L}^{\lambda(X)}})
& =BF^i_\lambda(X),\qquad i=0,1 \\ 
K_0({\cal F}_{{\frak L}^{\lambda(X)}})
& =\Delta^\lambda(X), 
\end{align}
where 
${\cal F}_{{\frak L}^{\lambda(X)}}$ is the AF-algebra defined 
by the fixed point algebra of 
${\cal O}_{{\frak L}^{\lambda(X)}}$
under the gauge action.
\end{thm} 
\begin{pf}
By Lemma 5.1 the  $\lambda$-synchronizing  $\lambda$-graph system  
${\cal L}^{\lambda(X)}$ 
of $X$ satisfies $\lambda$-condition (I).
Also, if $X$ is  
synchronized irreducible,
  then ${\cal L}^{\lambda(X)}$ is $\lambda$-irreducible,
 and by \cite[Theorem 3.9]{Ma2005b}  
${\cal O}_{{\frak L}^{\lambda(X)}}$ 
is simple.
The equalities 
(5.2),(5.3),(5.4) follow from 
(3.3),(3.4),(3.5) and (3.7),(3.8),(3.9)
(see also \cite{Ma1999a}, \cite{Ma2001b}, \cite{Ma2002a}).
\end{pf}

\begin{cor}
Let $X$ be a synchronizing subshift
satisfying synchronizing condition (I).
Then the  
 $C^*$-algebra
${\cal O}_{{\frak L}^{S(X)}}$
associated with the synchronizing $\lambda$-graph system 
${\frak L}^{S(X)}$ is a simple $C^*$-algebra  that
is isomorphic to the  $C^*$-algebra
${\cal O}_{{\frak L}^{\lambda(X)}}$
associated with the $\lambda$-synchronizing $\lambda$-graph system 
${\frak L}^{\lambda(X)}$ for $X$.
\end{cor}
\begin{pf}
As a synchronizing subshift is irreducible,
it is  synchronized irreducible
so that 
the  $\lambda$-synchronizing  $\lambda$-graph system  
${\cal L}^{\lambda(X)}$ 
of $X$ is $\lambda$-irreducible.
By Proposition 3.4, 
${\cal O}_{{\frak L}^{S(X)}}$
is isomorphic to 
${\cal O}_{{\frak L}^{\lambda(X)}}$
that is simple by the above theorem.
\end{pf}
One can prove  that 
the $C^*$-algebras 
${\cal O}_\beta, 1<\beta \in {\Bbb R},$ in \cite{KMW}
for the $\beta$-shifts 
$\Lambda_\beta$
are isomorphic to
the $C^*$-algebras 
${\cal O}_{{\frak L}^{\lambda(\Lambda_\beta)}}$
associated with the 
 the $\lambda$-synchronizing $\lambda$-graph systems for $\Lambda_\beta$.
 The $\lambda$-graph systems studied in the paper 
\cite{Ma9} are also the $\lambda$-synchronizing $\lambda$-graph systems
for the subshifts.



\end{document}